 \documentclass[final,5p,twocolumn]{elsarticle}

\usepackage{amssymb}
\usepackage{amsmath, amssymb}
\usepackage{amsfonts}
\usepackage{amsthm}
\usepackage{graphicx}
\usepackage{enumerate}
\usepackage[font = small]{caption}
\usepackage{ifthen}
\usepackage{multicol}
\usepackage{float}
\usepackage{fancyhdr}
\usepackage[usenames, dvipsnames]{color}
\usepackage[lofdepth,lotdepth]{subfig}
\graphicspath{{figures/}}
\usepackage{mathtools}

 \usepackage{tikz}
 \usetikzlibrary{decorations.pathmorphing,shadings}
\usepackage{tkz-graph}
\usetikzlibrary{arrows}

 \usepackage{cancel} 
 
\newtheorem{theorem}{Theorem}
\newtheorem{lemma}[theorem]{Lemma}
\newtheorem{proposition}[theorem]{Proposition}
\newtheorem{corollary}[theorem]{Corollary}
\newtheorem{result}[theorem]{Result}
\newtheorem{rem}{Remark}
\newtheorem{example}{Example}
\newtheorem{ass}{Assumption}
\newtheorem{definition}{Definition}

\newboolean{showcomments}
\setboolean{showcomments}{true}

\newcommand{\todo}[1]{  \ifthenelse{\boolean{showcomments}}
{\textcolor{ForestGreen}{TO DO:  #1}}{}}
\newcommand{\bassam}[1]{\ifthenelse{\boolean{showcomments}}
{\textcolor{Blue}{(Bassam says: #1)}}{}}
\newcommand{\emma}[1]{\ifthenelse{\boolean{showcomments}}
{\textcolor{VioletRed}{(Emma says: #1)}}{}}
\newcommand{\ifneeded}[1]{\ifthenelse{\boolean{showcomments}}
{\textcolor{Gray}{#1}}{}}

\newboolean{shownew}
\setboolean{shownew}{true}
\newcommand{\newtext}[1]{\ifthenelse{\boolean{shownew}}
{\textcolor{black}{#1}}{}}
\newboolean{edit}
\setboolean{edit}{true}
\newcommand{\edit}[1]{\ifthenelse{\boolean{edit}}
{{#1}}{}}

\newcommand*\wideestimates{\mathrel{\widehat{=}}}

\newcommand{\Vn}{V_N}

\newcommand{\hn}{$\mathcal{H}_2$ }

\newcommand{\Zld}{\mathbb{Z}_L^d}
\newcommand{\Znd}{\mathbb{Z}_L^d}
\newcommand{\xtr}{\hat{p}(\theta)}

\newcommand\TL{\rule{0pt}{3.5ex}}       
\newcommand\BL{\rule[-2ex]{0pt}{0pt}} 

\usepackage{array}
\newcolumntype{L}[1]{>{\raggedright\let\newline\\\arraybackslash\hspace{0pt}}m{#1}}
\newcolumntype{C}[1]{>{\centering\let\newline\\\arraybackslash\hspace{0pt}}m{#1}}
\newcolumntype{R}[1]{>{\raggedleft\let\newline\\\arraybackslash\hspace{0pt}}m{#1}}
\usepackage{subfig}

\usepackage{tikz}

 \biboptions{sort&compress}

\journal{Systems and Control Letters}

\begin{document}

\begin{frontmatter}



\title{Noise-Induced Limitations to the Scalability of Distributed Integral Control }


\author{Emma Tegling and Henrik Sandberg}

\address{School of Electrical Engineering and Computer Science, KTH Royal Institute of Technology, SE-100 44 Stockholm, Sweden \\ Corresponding author: E. Tegling, tegling@kth.se}

\begin{abstract}
We study performance limitations of distributed feedback control in large-scale networked dynamical systems. Specifically, we address the question of how the performance of distributed integral control is affected by measurement noise. We consider second-order consensus-like problems modeled over a toric lattice network, and study asymptotic scalings (in network size) of \hn performance metrics that quantify the variance of nodal state fluctuations.
While previous studies have shown that distributed integral control fundamentally improves these performance scalings compared to distributed proportional feedback control, our results show that an explicit inclusion of measurement noise leads to the opposite conclusion. The noise's impact on performance is shown to decrease with an increased inter-nodal alignment of the local integral states.
However, even though the controller can be tuned for acceptable performance for any given network size, performance will degrade as the network grows, limiting the \emph{scalability} of any such controller tuning. In particular, the requirement for inter-nodal alignment increases with network size. We show that this may in practice imply that very large and sparse networks will require any integral control to be centralized, rather than distributed. In this case, the best-achievable performance scaling, which is shown to be that of proportional feedback control, is retrieved.
\end{abstract}


\end{frontmatter}


\section{Introduction}
\label{sec:INTRO}
A central issue in the control of networked systems is to understand and quantify how the limited sensing, actuation and connectivity of a distributed controller structure affect global performance. A prototypical problem is that of distributed consensus, where the objective is to drive a network of agents to the same state, but where each agent only has access to limited and localized measurements. Natural questions arise as to how well a state of consensus can be upheld, for example, under external disturbances, and how this depends on the size of sensing neighborhoods and the topology of the controller. An understanding of these issues is key in achieving efficient and robust control performance in a wide range of applications, including vehicle platooning and formation control problems, wireless sensor networks and electric power systems.

In response to this issue, an ongoing research trend is to characterize fundamental limitations of distributed feedback control in terms of asymptotic bounds on various performance metrics~\cite{Bamieh2012,Lin2012,  Patterson2014, SiamiMotee2015,Grunberg2016,Barooah2009,herman2015nonzero}.
 In particular, the approach in~\cite{Bamieh2012} was to study distributed static state feedback controllers with locality constraints and derive scalings (in network size) of the best-achievable performance bounds. It was shown that a reasonable performance scaling in sparse networks requires that the local controllers have access to measurements of their own states with respect to a global reference frame, what is referred to as \emph{absolute feedback}. 
This observation motivated the work in~\cite{Tegling2017,Tegling2017b} where it was shown that for double-integrator networks, an absolute measurement of only \emph{one} of the two states (e.g. position \emph{or} velocity) can suffice. The addition of appropriately filtered distributed derivative or integral control can then namely alleviate the performance limitations that applied to static feedback. In this paper, we consider the same scenario, and focus on the distributed integral controller whose superior performance compared to distributed static feedback was shown in~\cite{Tegling2017}.

In line with standard intuition, integral control in networked dynamical systems is motivated by a desire to eliminate stationary control errors, and has been proposed in e.g.~\cite{Freeman2006, Andreasson2014, Seyboth2015,Lombana2015, Lombana2016}. In particular, it is important for frequency control in electric power networks, in order to reject disturbances and drive the system frequency to the desired setpoint (50 Hz or 60 Hz)~\cite{Andreasson2014ACC}. In that context, the integral action is referred to as \emph{secondary} frequency control. It is worth pointing out that while integral control can be implemented with various degrees of centralization, distributed approaches may be desirable (or the only feasible option) in many network applications. 

The question that has motivated the present work is to which extent the superior performance of distributed integral control compared to standard distributed static feedback is robust to measurement noise in the controller. 
The apparent reason for the improved performance is namely that integration of the absolute velocity measurements emulates absolute position feedback~\cite{Tegling2017}. Any noise and bias in the velocity measurements is prevented from causing destabilizing drifts in this position feedback by a \emph{distributed averaging filter} in the controller we consider.
 Yet, we show here that noisy measurements may still have a large impact on performance. 

Following the problem setup in~\cite{Bamieh2012, Tegling2017b} we consider networked systems modeled over toric lattices, where the local dynamics are of second order. We are concerned with the performance of these systems in terms of nodal variance measures that capture the notion of network \emph{coherence}, and evaluate how these measures scale asymptotically with the size of the network. An unfavorable scaling of nodal variance implies that performance will degrade as the network grows.  In such cases, the control law in question is limited in its \emph{scalability} to large networks. 

We show that while the performance of noiseless distributed integral control scales well, the addition of measurement noise gives rise to its own contribution to nodal variance with an unfavorable scaling. Even though this contribution, which is also proportional to the noise intensity, may be small in absolute terms for small networks it limits the overall scalability of the controller. In fact, it becomes even worse than with distributed static feedback. 

This paper extends the related work in~\cite{Flamme2018}, which treated this problem for electric power networks and with an alternative performance objective, and deepens the analysis. In particular, we here study the impact of the distributed averaging filter directly, \edit{and allow it to take on a different structure than the underlying feedback network.} We demonstrate that the inter-nodal alignment of integral states through this filter is important for performance.
 While this may seem intuitive, we show that the need for such alignment does not only increase with noise intensity, but more importantly, with the network size. In a 1-dimensional lattice, this increase is even faster than linear. 
This paper's main conclusion is therefore that scalable integral control in lattice networks can only be implemented in a centralized fashion, or must allow for \edit{a very high connectivity}.

The remainder of this paper is organized as follows. We introduce the problem formulation in~Section~\ref{sec:setup} and present the performance scalings with the various controllers in Section~\ref{sec:perf}. In Section~\ref{sec:h2density} we review the technical framework from~\cite{Tegling2017b} that is used to analyze the scalings in Section~\ref{sec:Aanalysis}. In particular, Section~\ref{sec:Aanalysis} treats the importance of the distributed averaging filter for the controller scalability. \edit{We present numerical examples in Section~\ref{sec:examples}} and conclude by a discussion of our findings in Section~\ref{sec:discussion}.

\section{Problem setup}
\label{sec:setup}
\subsection{Definitions}
Consider a network defined on the $d$-dimensional discrete torus $\Zld$. This is a lattice with a total of $N = L^d$ nodes and periodic boundary conditions. In the 1-dimensional case ($d = 1$), $\mathbb{Z}_L$ is simply the $N$ node ring graph.
 We will discuss \emph{scalings} of performance metrics with respect to the size of the network. The notation $\sim$ is used to denote scalings as follows:
\begin{equation}
\label{eq:scalingdef}
u(N) \sim v(N) ~~\Leftrightarrow~~ \underline{c}v(N) \le u(N) \le \bar{c}v(N),
\end{equation}
for any $N >0$, where the fixed constants $\underline{c},\bar{c}>0$ are independent of~$N$. When a scaling is said to hold \emph{asymptotically}, \eqref{eq:scalingdef} holds for all $N\ge \bar{N}$ for some~$\bar{N}$.

\subsection{System dynamics }
We treat a networked dynamical system where the local dynamics are of second order. This means that there are two $d-$dimensional states, $x_k$ and $v_k$, at each network site $k \in \Zld$. These states can be thought of as, respectively, the position and velocity deviations of the $k^\text{th}$ agent in a formation control problem, but may also capture, for example, phase and angular frequency in coupled oscillator networks (see Example~\ref{ex:powersys}). 
The system dynamics are modeled as follows (omitting the states' time-dependence in the notation):
\begin{equation}
\label{eq:system}
\begin{bmatrix}
\dot{x} \\ \dot{v}
\end{bmatrix} = \begin{bmatrix}
0 & I\\ F& G
\end{bmatrix} \begin{bmatrix}
x \\ v
\end{bmatrix} + \begin{bmatrix}
0 \\ I
\end{bmatrix}u+\begin{bmatrix}
0 \\ I
\end{bmatrix}w,
\end{equation}
where $u$ is a control input and $w$ \edit{models an uncorrelated disturbance entering at every network site}. The linear feedback operators $F$ and $G$ define convolutions of the states~$x$ and~$v$ with the function arrays~$f = \{f_k\}$ and~$g =  \{g_k\}$ over~$\Zld$, i.e., $(Fx)_k = \sum_{l \in \Zld} f_{k-l}x_l$.\footnote{Both the state and the function arrays are $d$-dimensional. The convolution is thus multi-dimensional and the multi-indices $k$ and $l$ are added as $k+l = (k_1,\ldots,k_d)+(l_1,\ldots,l_d) = (k_1+l_1,\ldots,k_d+l_d)$. To simplify the reading of this short letter, we will avoid the multi-index notation. A more detailed treatment of technicalities related to the states' dimensionality is found in~\cite{Tegling2017b}. } This structure implies that the state feedback is \emph{spatially invariant} with respect to $\Zld$. 
We refer to the system~\eqref{eq:system} as subject to \emph{static feedback} if the control input $u = 0$, since the feedback in this case is simply proportional to state deviations.\footnote{Alternatively, any control on the form $u = F^ux + G^uv$, where $F^u$ and $G^u$ satisfy Assumptions~\ref{ass:relative}--\ref{ass:decoupling} is possible. W.L.O.G. we can then assume $u=0$ and absorb $F^u$ and $G^u$ in~\eqref{eq:system}.}  
An example of the dynamics~\eqref{eq:system} is nearest-neighbor consensus for $d=1$:
\begin{multline}
\label{eq:exampleconsensus}
\ddot{x}_k \! = \!\dot{v}_k\! = \! f_+ \!\!\left( x_{k+1} \!-\!x_k\right) + f_- \!\!\left( x_{k-1} \!-\! x_k\right) + g_+\!\! \left( v_{k+1}\!-\!v_k \right) \\ + g_-\!\!\left( v_{k-1} \!-\! v_k\right) - f_o x_k - g_o v_k + u_k +w_k,
\end{multline}
where $f_+,f_-,f_o,g_+,g_-,g_o\ge 0$ are fixed gains.
We refer to terms like $\left( x_{k+1} \!-\!x_k\right)$ as \emph{relative feedback} and to terms like $-f_ox_k$ as \emph{absolute feedback}. Absolute feedback is well-known to be beneficial for control performance in networked dynamical systems, but the corresponding measurements are often not available (see e.g.~\cite{Barooah2007,Bamieh2012}). 
Here, we therefore make the following assumption on the system:
\begin{ass}[Relative position measurements]
\label{ass:relative}
Only relative measurements of the state $x$ are available, so the feedback can only involve differences between states of neighboring nodes. For the feedback operator~$F$, this implies that $\sum_{k \in \Zld} f_k = 0$ and in \eqref{eq:exampleconsensus} that $f_o = 0$. 
\end{ass}
That is, while each local controller has access to an absolute measurement of its (generalized) velocity, Assumption~\ref{ass:relative} implies that it cannot measure its position with respect to a global reference frame. Consider also the following example from electric power systems:
\vspace{-1mm}
\begin{example}[Frequency control in power networks]
\label{ex:powersys}
Synchronization in power networks is typically studied through a system of coupled swing equations. Under some simplifying assumptions, the linearized swing equation, also referred to as \emph{droop control}, 
can be written as: 
\begin{equation}
\label{eq:swingeq}
m \ddot{\theta}_k + d\dot{\theta}_k = -\sum_{j \in \mathcal{N}_k} b_{kj}(\theta_k - \theta_{j}) +P_{m,k} + u_k, 
\end{equation}
where $\theta_k$ is the phase angle and $\omega_k = \dot{\theta}_k $ the frequency deviation at node $k$, and $m$ and $d$ are, respectively, inertia and damping coefficients. The parameter $b_{kj} = b_{jk}$ is the susceptance of the $(k,j)^{\mathrm{th}}$ power line model, $\mathcal{N}_k$ is the neighbor set of node~$k$ and $P_{m,k}$ is a net power injection. \edit{Here, $u_k$ is called a secondary control input.}  
The dynamics~\eqref{eq:swingeq} can be cast as the system~\eqref{eq:system}, with $x \wideestimates \theta$, $v \wideestimates \omega $ and treating fluctuations in $P_{m,k}$ as the disturbance $w_k$. 
\end{example}

\vspace{-1mm}
We remark that the analysis here is not limited to nearest-neighbor feedback, but we assume that measurements are available from a neighborhood of width $2q$. As in~\cite{Bamieh2012,Tegling2017b}, we make the following additional assumptions:

\vspace{-1mm}
\begin{ass}[Locality]
\label{ass:locality}
All feedback operators use measurements from a local neighborhood of width $2q$, where the \emph{feedback window} $q$ is independent of $L$. For the feedback operator~$F$, 
this means that~$f_k = 0$~if~$|k|>q.$
\end{ass}

\vspace{-2mm}
\begin{ass}[Reflection symmetry]
\label{ass:symmetry}
The feedback interactions on $\Zld$ are symmetric around each site~$k$. 
For example in~\eqref{eq:exampleconsensus} this requires $f_+ = f_-$ and $g_+ = g_-$. 
\end{ass}

\vspace{-2mm}
\begin{ass}[Coordinate decoupling]
\label{ass:decoupling}
The feedback in each of the $d$ coordinate directions is decoupled from the components in the other coordinates. The array elements associated with all feedback operators are also isotropic. 
\end{ass}

\subsection{Distributed integral control}
\begin{figure}
\centering
\includegraphics[width = 0.32\textwidth]{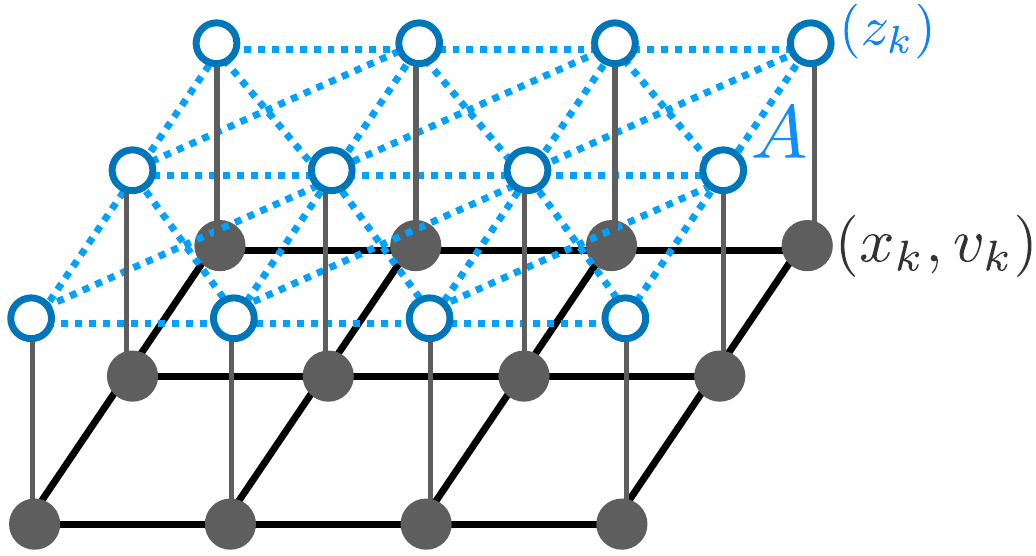}
\caption{ Example structure of the distributed integral controller. The inter-nodal alignment of integral states $z_k$ takes place over a communication network (dashed lines), while the state feedback interactions take place over the physical network (solid lines). {It is an example of a so-called {multiplex} network.}  }
\label{fig:networkfig}
\end{figure}

Consider the following control input to the system~\eqref{eq:system}:
\begin{equation}\label{eq:DAPI}
\begin{aligned}
u &= z \\ 
\dot{z} & = - c_ov^m + Az, 
\end{aligned}
\end{equation}
where $v^m$ is the velocity measured by the controller (for now, let $v^m = v$), $c_o>0$ is a fixed (integral) gain and $A$ is a feedback operator subject to the same assumptions as~$F$.  An example of the control law~\eqref{eq:DAPI} is:
\begin{equation}\label{eq:dapiexample}
\begin{aligned}
\dot{u}_k = \dot{z}_k = a_+(z_{k+1} - z_{k}) + a_-(z_{k-1} - z_k) - c_ov_k^m,
\end{aligned}
\end{equation}
where $a_+,a_->0$ are fixed gains. This controller integrates the absolute velocity measurements, but also aligns the integral state $z$ over the network through the consensus or \emph{distributed averaging filter} represented by the operator~$A$. The purpose of this alignment is to prevent drifts in the integral states~$z_k$ (due to noise or bias), which would otherwise destabilize the system~\cite{Andreasson2014ACC}. It is useful to think of the information exchange through~$A$ as taking place over a communication network layer, separate from the physical network. \edit{This layered structure results in what is sometimes referred to as a \emph{multiplex} network (see e.g.~\cite{Lombana2016}). }The setup is illustrated in Figure~\ref{fig:networkfig}.

This type of controller has been proposed in the context of power system frequency control in~\cite{Andreasson2014ACC, SimpsonPorco2013} for the elimination of stationary control errors that arise through standard droop control. Its key advantage is that it can be implemented in a distributed fashion. It is therefore more amenable to modern power networks with increasingly distributed generation than traditional, centralized secondary frequency control.


\begin{table*}[h]
\centering
\caption{Asymptotic performance scalings for the system~\eqref{eq:system} with (i) $u = 0$ (static feedback), (ii) $u$ as in~\eqref{eq:DAPI} with $v^m = v$ (distributed integral control, noiseless) and (iii) $u$ as in~\eqref{eq:DAPI} with $v^m = v +\varepsilon\eta$ (noisy distributed integral control). Scalings are up to a constant independent of network size $N$, algorithm parameter $\beta = \max \{||f||_\infty,||g||_\infty\}$ and relative noise intensity $\varepsilon$.
 }
\begin{tabular}{|l|l|l|}
\hline 
\rule[-1ex]{0pt}{2.5ex}  & \textbf{Local error }&\textbf{ Global error} \TL \BL \\ 
\hline 
\rule[-1ex]{0pt}{2.5ex} (i) Static feedback &$ \Vn \sim \frac{1}{\beta}$ for any $d$ & $ \Vn \sim \frac{1}{\beta}\begin{cases} N & d = 1\\  \log N & d = 2 \\  1& d \ge 3 \end{cases} $ \TL \BL\\ 
\hline 
\rule[-1ex]{0pt}{2.5ex} (ii) Distributed integral control (noiseless) & $ \Vn \sim \frac{1}{\beta}$ for any $d$ & $ \Vn \sim \frac{1}{\beta}$ for any $d$\TL \BL \\ 
\hline 
\rule[-1ex]{0pt}{2.5ex} (iii) Noisy distributed integral control &  $\Vn \sim   \frac{1}{\beta}+\frac{\varepsilon^2}{\beta}\begin{cases} N & d = 1\\  \mbox{log}N & d = 2 \\  1& d \ge 3 \end{cases}$ & 
$\Vn \sim \frac{1}{\beta}+ \frac{\varepsilon^2}{\beta}\begin{cases} N^3 & d = 1\\  N & d = 2 \\  N^{1/3}& d = 3 \\  \mbox{log}N & d = 4 \\ 1 & d\ge 5 \end{cases} $ \TL \BL\\ 
\hline 
\end{tabular} 
\label{tab:resultstab}
\end{table*}

\section{Performance of static feedback vs. distributed integral control}
\label{sec:perf}
We are concerned with the performance of the system~\eqref{eq:system}, and in particular, with how well the performance of a fixed control law scales as the network size $N \rightarrow \infty$. In line with related work \cite{Bamieh2012,Lin2012, Patterson2014, SiamiMotee2015,Grunberg2016, Tegling2017}, we characterize performance through the steady state variance of nodal state fluctuations, when the system is driven by a white noise disturbance input $w$  \edit{that is uncorrelated across the input channels, that is, the network sites}. 
For a system of size $N$, this variance can be characterized through the squared \hn norm from $w$ to a performance output $y$:
\begin{equation}
\label{eq:variancedef}
\mathbf{V}_N: = \sum_{k \in \Zld} \lim_{t\rightarrow \infty} \mathbb{E} \{y_k^T(t) y_k(t)\}.
\end{equation}
We consider the following performance measurements:

\vspace{-2.5mm}
\begin{definition}[Global error]
\begin{equation}
y_k = x_k - \frac{1}{N}\sum_{l\in \Zld} x_l
\label{eq:globalerror}
\end{equation}
This quantity measures the deviation of each state with respect to the network average and is therefore a measure of \textit{global} disorder.
\end{definition}

\vspace{-2.5mm}
\begin{definition}[Local error]
\begin{equation}
\label{eq:localerror}
y_k = x_k - x_{k-1}
\end{equation}
This quantity measures the deviation of each state with respect to its nearest neighbor and is therefore a measure of \textit{local} disorder.
\end{definition}
\vspace{-2mm}
Throughout this paper, we consider the \emph{per-site variance}, which is obtained by simply dividing the total \hn norm by the system size $N$. As the systems we consider are spatially invariant, the per-site variance is independent of the site $k$.
\begin{definition}[Per-site variance]
\label{def:persitevariance}
\begin{equation}
\Vn = \lim_{t\rightarrow \infty} \mathbb{E} \{y_k^T(t) y_k(t)\} =  \frac{\mathbf{V}_N}{N}.
\label{eq:persitevariance}
\end{equation}
\end{definition}
We are interested in the \emph{scaling} of the per-site variance~$\Vn$ with the system size~$N$ as it grows asymptotically. If~$\Vn$ scales slowly in $N$, we call the system more \emph{coherent} than one in which~$\Vn$ scales faster. 
It is only if the variance~$\Vn$ is \emph{bounded} in $N$ that we can say that a control law is \textit{scalable} to large networks.

The following results, of which (i) appeared in \cite[Corollary 3.2]{Bamieh2012} and (ii) follows from \cite[Corollary 1]{Tegling2017} are the main motivation for this work. 
\begin{result}[Performance scalings]
Consider the system~\eqref{eq:system} and let Assumptions~\ref{ass:relative}--\ref{ass:decoupling} hold. Assume that the velocity measurements are \emph{noiseless}, that is, $v^m = v$. Then, Table~\ref{tab:resultstab} lists the asymptotic scaling of the per-site variance~$\Vn$ with
\vspace{-1mm}
\begin{enumerate}[(i)]
\item Static feedback, i.e., where the secondary control input $u = 0$, and 

\vspace{-1mm}
\item Distributed integral control with $u$ given in~\eqref{eq:DAPI}. 
\end{enumerate}

\end{result} 
\begin{rem}
With respect to the local error, the distributed integral controller offers no improvement in terms of the scaling of the per-site variance compared to static feedback. In absolute terms, however, the variance is reduced (see Proposition~\ref{prop:h2expr} and note that $|\varphi(\theta)|>0$). 
\end{rem}



\subsection{Limitations due to noisy measurements}
Result~1 demonstrated that distributed integral control on the form~\eqref{eq:DAPI}, aside from its benefits in eliminating stationary control errors, can fundamentally improve performance in terms of the per-site variance of the global error. As discussed in~\cite{Tegling2017}, this improvement can be attributed to the fact that the integration of absolute velocity measurements can provide a substitute for the otherwise lacking absolute position feedback. It turns out, however, that this result is very sensitive to the accuracy of the absolute velocity measurements, and may change radically if they are subject to noise.

Here, let us therefore model additive measurement noise and let the velocity measurement in~\eqref{eq:DAPI} be
\[v^m = v + \eta,\]
where the vector $\eta$ contains uncorrelated white noise with the relative intensity~$\varepsilon$ defined through $\mathbb{E}\{\eta(\tau)\eta^T(t)\} = \varepsilon\mathbb{E}\{w(\tau)w^T(t)\}$. 
Inserting into~\eqref{eq:system} gives: 

\begin{equation}
\begin{aligned}
\begin{bmatrix}
\dot{z} \\ \dot{x} \\ \dot{v}
\end{bmatrix} = \begin{bmatrix}
A & 0 & -c_oI \\ 0& 0 & I\\ I & F& G
\end{bmatrix} \begin{bmatrix}
z\\ x \\ v
\end{bmatrix} + \begin{bmatrix}
0& -c_o\varepsilon I\\ 0 &0 \\ I & 0
\end{bmatrix} \bar{w},
\end{aligned}
\label{eq:noisydapivector}
\end{equation}
where $\bar{w} \in \mathbb{R}^{2N}$ is a vector of uncorrelated white noise.
Evaluating local and global performance scalings for this system leads to the following result. 

\vspace{-1.5 mm}
\begin{result}[Performance scalings with noise]
\label{prop:noisyscaling}
Consider the system~\eqref{eq:noisydapivector} and let Assumptions~\ref{ass:relative}--\ref{ass:decoupling} hold. Then, row~(iii) of Table~\ref{tab:resultstab} lists the asymptotic scaling of the per site variance $\Vn$. 
\end{result}
\vspace{-3mm}
\begin{proof}
Follows from \edit{the upcoming} Proposition~\ref{prop:h2expr} and Corollaries~\ref{cor:pw}--\ref{cor:peta}.
\end{proof}
\vspace{-1.5 mm}

Result~\ref{prop:noisyscaling} reveals that the measurement noise~$\eta$ leads to an unfavorable scaling of both local and global error variance -- even worse than with static feedback. This may not be an issue for small networks, as the variance is scaled by the factor $\varepsilon^2$, which can be very small (recall, $\varepsilon$ represents the intensity of the measurement noise~$\eta$ relative to the process disturbance $w$). However, performance will deteriorate as the network size grows, thus limiting the scalability of distributed integral control.

\vspace{-1.5mm}
\begin{rem}
Here, we have assumed that the velocity enters without noise in the system dynamics~\eqref{eq:system}. It may also be reasonable to model the same noise there, so that $\dot{v} = Fx + G(v+\eta) + u +w$. This can, however, be shown not to affect the qualitative system behavior discussed here~\edit{\cite{TeglingThesis}}. 
\end{rem}
\vspace{-1.5 mm}



\section{The \hn norm density and asymptotic performance scalings}
\label{sec:h2density}
We now review the technical results that were used to derive Table~1, and which will be needed to further analyze the impact of control design on performance. These results can all be found in~\cite{Tegling2017b} along with a more detailed discussion. 

\subsection{Diagonalization using Fourier transforms}
The systems considered in this paper can all be block-diagonalized by the spatial discrete Fourier Transform~(DFT). For a feedback operator $F$ with associated function array $f:~\Zld \rightarrow \mathbb{R}$, this is defined as
\( \hat{f}_n := \sum_{k \in \Zld} f_ke^{-i \frac{2\pi}{L}n\cdot k},\)
where $n = (n_1,\ldots,n_d)$ is a wavenumber.
All feedback operators considered herein are local by Assumption~\ref{ass:locality}. They can therefore be unambiguously re-defined onto the infinite lattice $\mathbb{Z}^d$ by adding zero entries wherever $|k|>q$. 
The $Z$-transform can then be taken as
\(\hat{f}(\theta) := \sum_{k \in \mathbb{Z}^d} f_ke^{-i \theta \cdot k},\)
where $\theta = (\theta_1,\ldots,\theta_d) \in [-\pi,\pi]^d$ is a spatial frequency. 

It is now easy to see that the DFT is sub-samples of the $Z$-transform at each wavenumber:
\begin{equation}
\hat{f}_n = \hat{f}\left(\theta = \frac{2\pi}{L}n\right),~~n\in\Zld.
\end{equation}  
We refer to $\hat{f}_n$ and $\hat{f}(\theta)$ as (generalized) \textit{Fourier symbols}.
For the general state-space system 
\begin{equation}
\begin{aligned}
\label{eq:generalss}
\dot{\psi} & = \mathcal{A}\psi + \mathcal{B}w\\
y & = \mathcal{C}\psi
\end{aligned}
\end{equation}
we can obtain the matrix-valued DFTs $\hat{\mathcal{A}}_n,~ \hat{\mathcal{B}}_n,~ \hat{\mathcal{C}}_n$, which are subsamples of the $Z$-transforms $ \hat{\mathcal{A}}(\theta),~ \hat{\mathcal{B}}(\theta),~ \hat{\mathcal{C}}(\theta)$. The eigenvalues of $\mathcal{A}$ are then simply all eigenvalues of $\hat{\mathcal{A}}(\theta)$ as $\theta = \frac{2\pi }{L}n,$ $n\in\Zld$.
\begin{example}
\label{ex:diagonalization}
For the system~\eqref{eq:system} with static feedback ($u = 0$) we have 
\[\mathcal{A}(\theta) = \begin{bmatrix}
0 & 1\\ \hat{f}(\theta) & \hat{g}(\theta)
\end{bmatrix},~~\mathcal{B}(\theta) = \begin{bmatrix}
0 \\ 1
\end{bmatrix}\]
where $\hat{f}(\theta) = -\sum_{k \in \mathbb{Z}^d}f_k(1 - \cos(\theta \cdot k))$ and $\hat{g}(\theta) = -g_o -\sum_{k \in \mathbb{Z}^d}g_k(1 - \cos(\theta \cdot k))$. 
\end{example}

For the output measurement~\eqref{eq:globalerror} we have that $\left. \hat{\mathcal{C}}(\theta) = \begin{bmatrix}1 & 0\end{bmatrix} \right.$ for $\theta \neq 0$. For the local error it holds $\hat{\mathcal{C}}^*(\theta)\hat{\mathcal{C}}(\theta) = \begin{bmatrix}\hat{l}(\theta) & 0 \\ 0& 0\end{bmatrix}$ with $\hat{l}(\theta) = 2(1 \!-\!\cos\theta)$. 
In both cases, $\hat{\mathcal{C}}(\theta =0) = 0$ and the subsystem $(\hat{\mathcal{A}}(0),\hat{\mathcal{C}}(0))$ is therefore unobservable. 

\subsection{\hn norm evaluation}
Provided that $\hat{\mathcal{A}}(\theta)$ is Hurwitz for all $\theta \neq 0$, the per-site variance $\Vn$ from~\eqref{eq:persitevariance} can be evaluated as 
\begin{equation} \label{eq:samplingsum}
\Vn = \frac{1}{N} \sum_{\mathclap{\substack{ \theta = \frac{2\pi}{L}n \\ n \in \Znd \backslash \{0\}} } } \mathrm{tr} \left( \hat{\mathcal{B}}^* (\theta) \hat{P}(\theta) \hat{\mathcal{B}}(\theta) \right), 
\end{equation}
where the observability Gramian $\hat{P}(\theta)$ at each $\theta \neq 0$ can be obtained by solving the Lyapunov equation
\begin{equation}
\label{eq:lyap}
\hat{\mathcal{A}}^*(\theta)\hat{P}(\theta)+\hat{P}(\theta)  \hat{\mathcal{A}}(\theta) = -\hat{\mathcal{C}}^*(\theta)\hat{\mathcal{C}}(\theta).
\end{equation}
The summand in~\eqref{eq:samplingsum} captures the distribution of the per-site variance $\Vn$ over the spatial frequency~$\theta$. We will therefore refer to it as the (per-site) \hn norm density:
\begin{definition}[Per-site \hn norm density]
\[\xtr =  \mathrm{tr} \left( \hat{\mathcal{B}}^* (\theta) \hat{P}(\theta) \hat{\mathcal{B}}(\theta) \right).\]
\end{definition}

\subsection{Bounds on asymptotic scalings}
The behavior of the \hn norm density determines the scaling of the per-site variance with the network size $N$. In particular, if $\xtr$ is uniformly bounded for $\theta \in [-\pi,\pi]^d$, then $\Vn$ is bounded in $N$. However, $\xtr$ has a singularity at $\theta = 0$ if $\hat{\mathcal{A}}(0)$ is non-Hurwitz. For example, $\hat{\mathcal{A}}(0)$ is non-Hurwitz with static feedback as a consequence of Assumption~\ref{ass:relative}. While the point at $\theta = 0$ is excluded from the sum in~\eqref{eq:samplingsum}, the singularity causes an unfavorable scaling of $\Vn$. Consider the following Lemma:
\begin{lemma}\emph{\cite[Lemma 4.2]{Tegling2017b}}
\label{lem:scalinglemma}
Assume the \hn norm density is such that
\begin{equation}
\label{eq:densityscaling}
\xtr \sim \frac{1}{\beta^{p}}\cdot\frac{1}{|\theta|^r}
\end{equation}
for $\theta \in [-\pi,\pi]^d$, where $\beta$ is an algorithm parameter, $p$ and $r$ are constants, and $|\cdot|$ denotes the Euclidean norm. Then, the per-site variance $\Vn$ scales asymptotically as
\begin{equation}
\label{eq:scalinggeneral}
\Vn \sim \frac{1}{\beta^{p}} \begin{cases}
 L^{r-d}  & ~~\mathrm{if}~ d< r \\
\log L  & ~~\mathrm{if}~ d= r \\
1 & ~~\mathrm{if}~ d >r.
 \end{cases}
\end{equation}
\end{lemma}
The systems considered in this paper all have \hn norm densities that can be written as in~\eqref{eq:densityscaling} with $r\in \{0,2,4\}$.
 To show this, the following Lemma is needed:
\begin{lemma}\emph{\cite[Lemma 6.3]{Tegling2017b}}
\label{lem:scalinglem}
For any admissible $F$ that satisfies Assumptions~\ref{ass:relative}--\ref{ass:decoupling}, it holds
\begin{equation}
\hat{f}(\theta) \sim -\beta|\theta|^2,
\end{equation}
where $\beta = ||f||_\infty$.  For~$G$, which contains absolute feedback, it holds $\hat{g}(\theta) \sim -g_o.$
\end{lemma}
\vspace{-0.7mm}
Any feedback operators (such as $A$) subject to the same assumptions as $F$ or $G$ have the same behavior. Therefore, the Fourier symbol for the local error measurement satisfies~$\hat{l}(\theta) \sim |\theta|^2$.
\begin{example}
\label{ex:scalingex}
Consider the system from Example~\ref{ex:diagonalization}. By solving the Lyapunov equation~\eqref{eq:lyap} with the outputs~\eqref{eq:globalerror} and \eqref{eq:localerror}, we obtain the \hn norm densities as
\begin{equation}
\label{eq:statich2density}
 \hat{p}^\mathrm{global}(\theta) = \frac{1}{2\hat{f}(\theta)\hat{g}(\theta)},~~\hat{p}^\mathrm{local}(\theta) = \frac{\hat{l}(\theta)}{2\hat{f}(\theta)\hat{g}(\theta)} 
\end{equation}
Lemma~\ref{lem:scalinglem} reveals that $\hat{p}^\mathrm{global}(\theta) \sim \frac{1}{\beta|\theta|^2}$ and $\hat{p}^\mathrm{local}(\theta) \sim \frac{1}{\beta}$. The scalings in Table~\ref{tab:resultstab} for static feedback then follow from Lemma~\ref{lem:scalinglemma}. 
\end{example}

\section{Improving the scalability of integral control}
\label{sec:Aanalysis}
\begin{figure}
\centering
	\includegraphics[width = 0.46\textwidth]{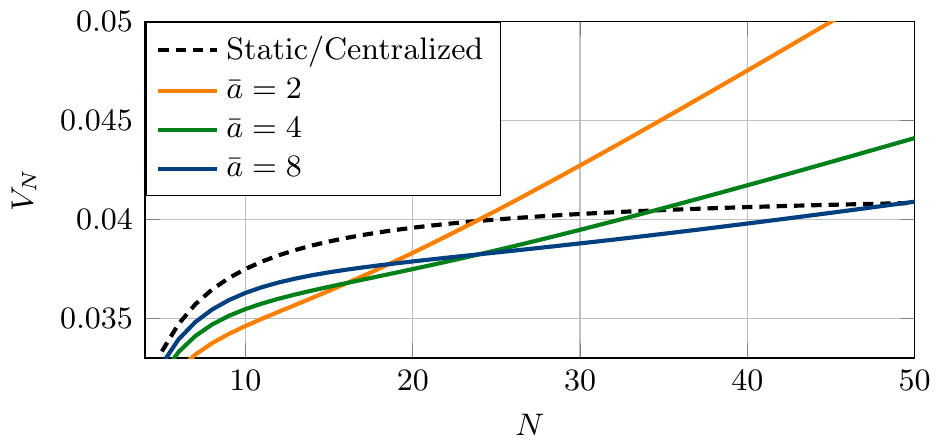}
\caption{Scaling of local error variance with static feedback vs. noisy distributed integral control in 1D lattice. For a given network size~$N$, it is possible to set the gain $\bar{a}$ in the distributed averaging filter $A$ so that the integral controller performs better than static feedback. Yet, no such controller scales well in~$N$. A centralized integral controller on the form~\eqref{eq:capi}, which corresponds to $\bar{a} \rightarrow \infty$, will however have the same performance as static feedback for any~$N$.  }
\label{fig:scalingplot}
\end{figure}

Let us now consider a situation where the system~\eqref{eq:system} is fixed and the design of the distributed integral controller~\eqref{eq:DAPI} for performance is of interest. In Proposition~\ref{prop:h2expr} below, we show that the error variance consists of two terms due to, respectively, disturbances and measurement noise.  For any given system of a fixed network size, it is possible to trade off these terms and to optimize the control design, as was the focus in~\cite{Flamme2018}. 
However, the unfavorable scaling of the error variance due to measurement noise sets fundamental limitations to the scalability of any such control design to large networks. A numerical example showcasing this issue is shown in Figure~\ref{fig:scalingplot}. The objective of this work, rather than to solve a performance optimization problem for a given system, is to point to the underlying limitations. 

To enable this analysis, consider the following \hn norm density expressions
\begin{proposition}
\label{prop:h2expr}
The \hn norm density of the system~\eqref{eq:noisydapivector} with respect to the \emph{global} error measurement~\eqref{eq:globalerror} is:
\begin{equation}
\label{eq:h2global}
\xtr \!=\! \underbrace{\frac{1}{2\hat{f}(\theta)} \! \cdot \! \frac{1}{ \varphi(\theta)+ \hat{g}(\theta) }}_{\hat{p}^w(\theta)} + \underbrace{ \frac{\varepsilon^2}{\hat{a}(\theta)}\! \cdot \!\frac{c_o}{2\hat{f}(\theta)} \! \cdot \!\frac{1}{1 + \hat{g}(\theta)\varphi^{-1}(\theta)} }_{\hat{p}^\eta(\theta)}.
\end{equation}
The \hn norm density with respect to the \emph{local} error measurement~\eqref{eq:localerror} is:
\begin{equation}
\label{eq:h2local}
\xtr \!=\! \underbrace{\frac{\hat{l}(\theta)}{2\hat{f}(\theta)} \! \cdot \! \frac{1}{  \varphi(\theta)+ \hat{g}(\theta) }}_{\hat{p}^w(\theta)} + \underbrace{ \frac{\varepsilon^2}{\hat{a}(\theta)}\! \cdot \!\frac{c_o\hat{l}(\theta)}{2\hat{f}(\theta)} \! \cdot \!\frac{1}{1+ \hat{g}(\theta)\varphi^{-1}(\theta) } }_{\hat{p}^\eta(\theta)}
\end{equation}
where 
\[ \varphi(\theta) = \frac{c_o (\hat{a}(\theta)  + \hat{g}(\theta) )}{\hat{a}^2(\theta)  +\hat{g}(\theta) \hat{a}(\theta)  -\hat{f}(\theta) }.\]
Here, $\hat{p}^w(\theta)$ corresponds to the \hn norm density of the system with noiseless distributed integral control and $\hat{p}^\eta(\theta)$ represents the contribution from the measurement noise. 
\end{proposition} 
\vspace{-4.5mm}\begin{proof}
The result follows from diagonalizing the system~\eqref{eq:noisydapivector} through Fourier transforms in line with Example~\ref{ex:diagonalization}, and then solving the corresponding Lyapunov equation~\eqref{eq:lyap}. The contributions from the disturbance inputs~$w$ and~$\eta$ can be separated since they are uncorrelated. 
\end{proof}
\vspace{-1mm}
The following corollaries lead to the results in Table~\ref{tab:resultstab}.
\vspace{-1mm}
\begin{corollary}
\label{cor:pw}
It holds that $\hat{p}^w(\theta)$ is uniformly bounded with respect to both global and local error, that is, $r = 0$ in Lemma~\ref{lem:scalinglemma}. 
\end{corollary}
\vspace{-4.5mm}
\begin{proof}
Substituting the scalings from Lemma~\ref{lem:scalinglem} into the expressions in Proposition~\ref{prop:h2expr} reveals that $\varphi(\theta) \sim {1}/{\beta|\theta|^2}$. The product $\hat{f}(\theta)\varphi(\theta)$ is thus bounded away from zero and the result follows. 
\end{proof}
\begin{corollary}
\label{cor:peta}
It holds $\hat{p}^\eta(\theta)\sim {\varepsilon^2}/{\beta |\theta|^4}$, or $r = 4$ (global error), and $\hat{p}^\eta(\theta)\sim {\varepsilon^2}/{\beta |\theta|^2}$, or $r = 2$ (local error). 
\end{corollary}
\vspace{-4.5mm}
\begin{proof}
Lemma~\ref{lem:scalinglem} gives that $\varphi^{-1}(\theta) \sim \beta |\theta|^2$. Since $\hat{f}(\theta) \sim -\beta |\theta|^2$, $\hat{a}(\theta) \sim \bar{a} |\theta|^2$, the product $\hat{f}(\theta)\hat{a}(\theta) \sim -\beta |\theta|^4$ and the result follows.  
\end{proof}


\subsection{From distributed to centralized integral control}
\edit{Under the given assumptions, the performance scalings in Table~\ref{tab:resultstab} hold with any design of the integral controller~\eqref{eq:DAPI}.  That is, for any fixed, finite gain $c_o$ and any operator $A$ with fixed, finite gains and subject to a locality constraint.  We now inquire whether better scalings can be achieved if these assumptions were relaxed. And if so, how must the controller be adjusted?}
The following conclusions can be drawn from Proposition~\ref{prop:h2expr}:
\vspace{-0.9mm}
\begin{enumerate}[a.]
\item It is not possible to set $\hat{a} = 0$ as in that case, $\hat{p}^\eta(\theta) = \infty$. 

\vspace{-0.9mm}
\item If $\hat{a}(\theta) \rightarrow \infty$, or $c_o \rightarrow 0$, then $\hat{p}^\eta(\theta) \rightarrow 0$, that is, the noise contribution to the variance vanishes. \\
At the same time, $\varphi(\theta)\rightarrow 0$ and $\hat{p}^w(\theta)$ becomes as with static feedback (compare \eqref{eq:h2global}--\eqref{eq:h2local} to~\eqref{eq:statich2density}).

\vspace{-0.9mm}
\item If $\hat{a}(\theta)$ is bounded away from zero, then $\hat{p}^\eta(\theta) \sim {\varepsilon^2}/{\beta|\theta|^2}$ (global error) and $\hat{p}^\eta(\theta) \sim {\varepsilon^2}/{\beta}$ (local error). \\
At the same time, $\varphi(\theta)$ becomes uniformly bounded and $\hat{p}^w(\theta) \sim {1}/{\beta|\theta|^2}$ (global error), while $\hat{p}^w(\theta) \sim {1}/{\beta}$ (local error), that is, the same as with static feedback.
\end{enumerate}
\vspace{-1.5mm}
Using these observations, the following result is derived:
\vspace{-1.5mm}
\begin{proposition}
\label{cor:bestscaling}
The best-achievable performance scaling for the noisy integral controlled system~\eqref{eq:noisydapivector} is that of distributed static feedback in Table~\ref{tab:resultstab}. 
\end{proposition}
\vspace{-5mm}
\begin{proof}
First, note that for any fixed $c_o>0$ and $\hat{a}(\theta)$, the scalings in Table~\ref{tab:resultstab} hold. 
For a better performance scaling, the behavior of $\hat{p}^\eta(\theta)$ in $\theta$ must change for the better ($r$ in Lemma~\ref{lem:scalinglemma} must decrease). This can only happen if $c_o \rightarrow 0$, $\hat{a}(\theta)\rightarrow \infty$ or if $\hat{a}(\theta)$ becomes bounded away from zero. As $\hat{p}^w(\theta)$ and $\hat{p}^\eta(\theta)$ have inverse dependencies on the function~$\varphi(\theta)$ in which both $c_o$ and $\hat{a}(\theta)$ appear, this will lead to cases b and c above. 
\end{proof}
\vspace{-3mm}
This means that the system can never have bounded variance in terms of the global error measurement~\eqref{eq:globalerror}. However, a bounded variance and thus scalability in terms of the local error~\eqref{eq:localerror} can be achieved \edit{by a re-tuning of the controller. 
Analyzing the cases} b and c above shows that the best-achievable performance scaling can be retrieved in three ways, which we discuss next.
\vspace{-2mm}
\subsubsection{Decreasing the integral gain $c_o$}
Decreasing the gain~$c_o$ reduces the impact of the measurement noise~$\eta$.  
To counteract the unfavorable scaling of~$\hat{p}^\eta(\theta)$, it must be ensured that $c_o/\hat{a}(\theta)$ is uniformly bounded in $\theta$. Since by Lemma 4, $\hat{a}(\theta) \sim -\bar{a}|\theta|^2$, this requires $c_o \sim \min|\theta|^2$.

The smallest wavenumber that contributes to the error variance in~\eqref{eq:samplingsum} corresponds to $\theta_{\min} = {2\pi}/{L}$. This implies that $c_o$ must be \emph{decreased as $1/L^2 $}. 
As the network grows, this implies  $c_o \rightarrow 0$ and the integral action is eliminated. In this case, the control input~$u$ is simply not used.

\vspace{-2mm}
\subsubsection{Increasing the distributed averaging gain}
For a fixed~$c_o$, the distributed averaging gain can be increased so that $\hat{a}(\theta)$ becomes bounded away from zero even as $L$ increases. Recall that $\hat{a}(\theta) \sim -\bar{a}|\theta|^2$ where $\bar{a} = ||a||_\infty$. This need not approach zero if $\bar{a} \sim 1/|\theta|^2$.

Again, $\theta_{\min} = {2\pi}/{L}$, meaning that $\bar{a}$ must be \emph{increased as $L^2$}. This implies that we must require~$\bar{a} \rightarrow \infty$ when the lattice size~$L$ grows.

While an infinite gain in distributed averaging is not feasible in practice, the same result can be realized as \emph{centralized} averaging integral control. Here, a central controller has instantaneous access to the integral states at all nodes. The control signal~$u_k$ is then the same for all $k \in \Zld$:
\begin{equation}
\label{eq:capi}
\begin{aligned}
u_k &= z;\\
\dot{z} &= \frac{1}{N} \sum_{k \in \Zld} v_k^m.
\end{aligned}
\end{equation}
It is not difficult to show that this controller has the same performance with respect to the errors~\eqref{eq:globalerror} and \eqref{eq:localerror} as static feedback.

\vspace{-2mm}
\subsubsection{Increasing communication network connectivity}
\label{sec:increasingcomms}
By relaxing Assumption~\ref{ass:locality} \edit{of locality }for~$A$, we can also bound $\hat{a}(\theta)$ away from zero. 
\edit{Let $q_A:=\max_{a_k \neq 0}|k|$ define the width of the feedback window in the communication network and consider the following lemma:}
\edit{ 
\vspace{-1.5mm}
\begin{lemma} 
\label{lem:q}
If $q_A \sim L^{2/3}$, then $|\hat{a}\left(\theta = {2\pi}/{L} \right)| \ge \delta$ for any~$L$, where $\delta$ is a positive constant. 
\end{lemma}
\vspace{-3.5mm}
\begin{proof}
See Appendix.
\end{proof}
\vspace{-1.5mm} This means that if the connectivity of the communication network is allowed to scale as $q_A\sim L^{2/3} = N^{2/3d}$, then $\hat{a}(\theta)$ stays bounded away from zero as $\theta \rightarrow 0$.
}
%
%


Allowing $q_A$ to increase with the lattice size~$L$ \edit{implies that new connections must be established as the network grows, and that a very high connectivity is required in large networks. }
This is practically challenging for large networks, and a centralized approach may be preferable.

\subsection{Implications for distributed integral control}
This section has shown that the distributed averaging filter~$A$ in the controller~\eqref{eq:DAPI} is important for performance. 
Recall that the role of the filter~$A$ is to align the controllers' integral states $z_k$ across the network, in order to gain robustness to measurement noise and bias. 
Previous results reported in~\cite{Wu2016,Andreasson2017,Tegling2017} have indicated that ``little'' inter-nodal alignment (i.e., small gains $\bar{a}$ and few interconnections in the communication network) is optimal for performance in the absence of measurement noise. 
It is intuitively clear, that the inter-nodal alignment through~$A$ becomes increasingly important if measurement noise is considered explicitly. 

Our results, however, reveal that it is not enough to scale the distributed averaging gain $\bar{a}$ 
with the noise intensity, here parameterized through~$\varepsilon$. Perhaps surprisingly, the need for inter-nodal alignment instead grows with the network size. It is required that the distributed averaging gain $\bar{a}\sim L^2 = N^{2/d}$, \edit{which for large, sparse networks in principle requires centralized integral control. Alternatively, the feedback window must be scaled so that $q_A \sim L^{2/3} = N^{2/3d}$. 
We demonstrate some of these results in the next section. }

\section{{Numerical examples}}
\label{sec:examples}
\edit{We present two numerical examples to illustrate the implications of this paper's main results in applications.}

\begin{figure*}
  \centering
  \includegraphics[width  =1\textwidth]{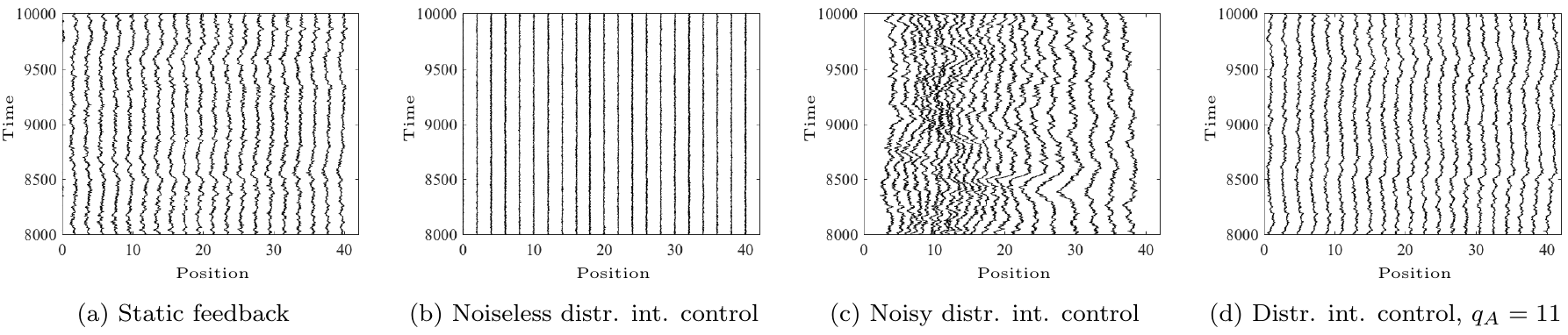}
      \caption{
{{Subset of position trajectories in a 100 vehicle platoon regulated according to~\eqref{eq:exampleconsensus}. Each vehicle is subjected to an independent random forcing at each time step in the simulation (the transients are not displayed). Under noiseless distributed integral control (b), the formation remains rigid with well-regulated spacings under such forcings -- it is coherent. The coherence is lost if the velocity measurements from the speedometers are noisy (c). Increasing the connectivity of the controller's communication network from $q_A=1$ to $q_A = 11$ as in (d) approximately retrieves the performance of static feedback (a).   }}}
\label{fig:simulation}
\end{figure*}

\subsection{{Lack of coherence in vehicle platoons}}
\edit{The performance limitations discussed in this paper are most severe for systems built over 1-dimensional lattice networks. Vehicle platooning problems, where vehicles are driven in a string formation, can therefore be expected to be fundamentally challenging if they are large. Here, we consider a platooning problem where each vehicle controls its velocity according to~\eqref{eq:exampleconsensus}. The position $x_k$ and velocity~$v_k$ represent deviations from the $k^{\text{th}}$ vehicle's desired trajectory $\bar{x}_k(t) = \bar{v}t+k\Delta$, where $\bar{v}$ is the cruising velocity and $\Delta$ is a desired inter-vehicle spacing. The disturbances~$w_k$ capture random forcings on each vehicle. }

\edit{In Figure~\ref{fig:simulation}, we have simulated an $N =100$
vehicle platoon controlled with (a) static feedback, (b) \emph{Noiseless} distributed integral control, and (c) \emph{Noisy} distributed integral control. 
As predicted by Table~\ref{tab:resultstab}, the formation is fully coherent under distributed integral control, provided the measurements are noiseless. If they are not, it causes a lack of coherence in the platoon that manifests as wave-like motions across the length of the platoon (see also~\cite{Bamieh2012}). The same behavior is visible, but to a lesser degree, under static feedback.}

\edit{In the simulation (d), we have allowed the communication network connectivity to increase in order to mitigate the impact of the measurement noise. Here, the size of the feedback window is $2q_A =22$. As predicted by Section~\ref{sec:increasingcomms}, this gives a performance similar to static feedback. This performance could be maintained also if the platoon grows by letting $q_A \sim N^{2/3}$.    }

\edit{Note that while the theoretical results in this paper are derived for spatially invariant systems (lattices with periodic boundary conditions), the simulated vehicle platoon has no connection between the first and last vehicles. It is also heterogeneous in the sense that the coupling strengths $f_{+},f_{-},g_{+},g_{-}$ between all vehicles were assigned random values in the interval $[0.5,1.5]$.  The simulation shows that our results seem robust to such violations of the spatial invariance assumption.   }
   

\begin{figure*}
  \centering
\includegraphics[width = \textwidth]{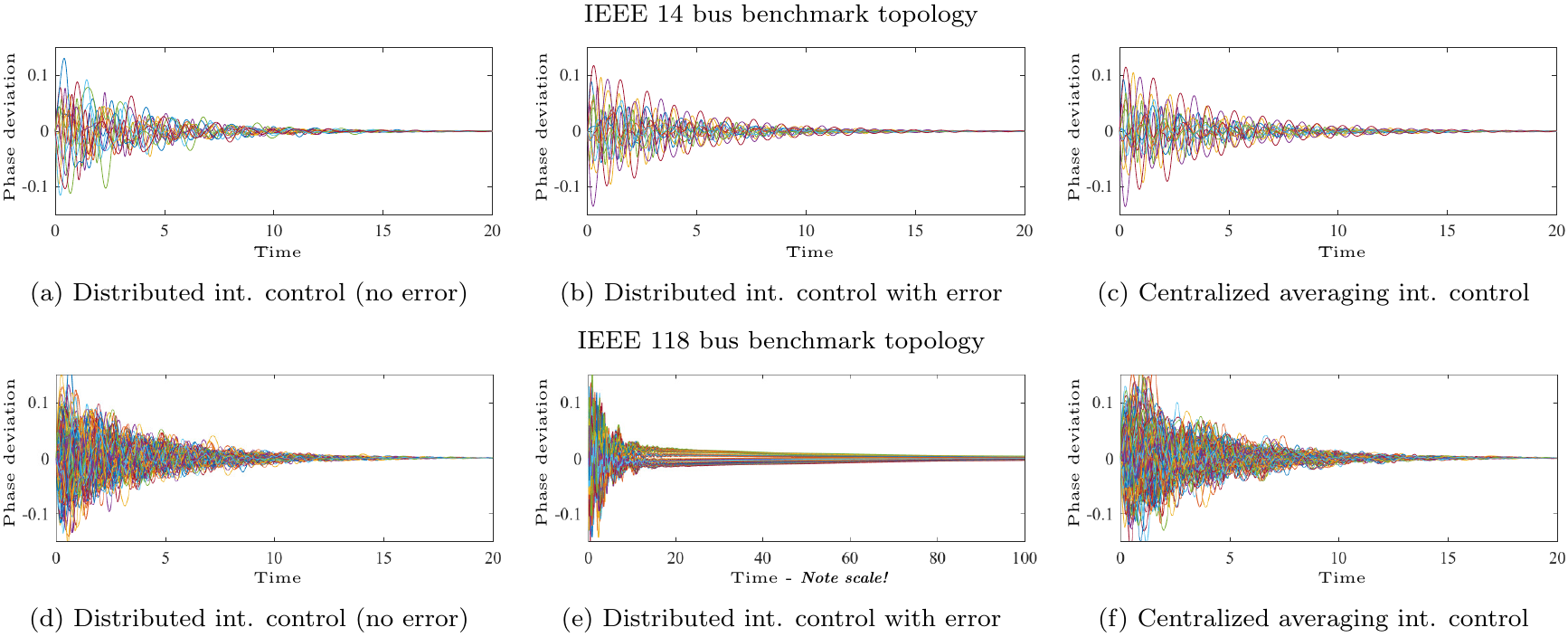}
\caption{Simulation of synchronization transients in a power network where a generator has been placed at each node of the IEEE 14 and 118 bus network topologies. The plots display deviations of generator phase angles from their respective setpoints after an initial perturbation in the net power loads and, in (b),(c),(e), and (f), an error in each integral controller’s frequency measurement. The measurement error limits the scalability of the distributed integral controller. It therefore has a significant impact on the transient performance of the 118 bus network, see~(e) and note the scale. Reverting to centralized averaging integral control as in~\eqref{eq:capi} improves performance, see (f).}
\label{fig:simulation2}
\end{figure*}

\subsection{{Transient performance in power networks}}
\edit{The distributed integral controller we consider has, as previously mentioned, been proposed for secondary frequency control in power networks.  Our results imply that this controller's scalability is limited if its frequency measurements are imperfect. As an example, consider the droop controller from Example~\ref{ex:powersys}. In this context, the distributed integral controller can be written
\begin{align*}
u_k & = z_k\\ 
\dot{z}_k &= -c_o\omega^m_k - \sum_{j \in \mathcal{N}^A_k}a_{ij}(z_{k}- z_j),
\end{align*}
where $\mathcal{N}^A_k$ denotes the neighbor set of node $k$ in the communication network.  }

\edit{Figure~\ref{fig:simulation2} displays synchronization trajectories for a hypothetical set of identical power generators ($m = \frac{20}{2\pi\cdot 60}$ and $d = \frac{10}{2\pi\cdot 60}$) placed at each node of, respectively, the IEEE 14 bus and 118 bus benchmark topologies~\cite{TestCases}, after an initial perturbation in the net power load. We compare the transient when the controllers have perfect frequency measurements to the case when the measurements have an initial error.   }

\edit{The scenario in the simulation corresponds to an alternative interpretation of the (squared) \hn norm, namely as the time integral of the system response \(\int_0^\infty \mathbb{E} \{y^T(t) y(t)\} \) when the initial condition $\psi_0$ is a random variable with $\mathbb{E}\{\psi_0 \psi_0^T\} = \mathcal{B}\mathcal{B}^T$ and the input is zero (see e.g.~\cite{Tegling2014}). Recall that the variance $V_N$ from Definition~\ref{def:persitevariance} represents a normalized \hn norm. Therefore, the scalings in Table~\ref{tab:resultstab} also offer insights to the simulated scenario.}

\edit{The simulation shows that the impact of frequency measurement errors on the transient performance is many times greater in the larger 118 bus network than in the small 14 bus network (compare, respectively, (a) to (b) and (d) to (e) and note the time scale). Here, the topology of the communication network is identical to the underlying power network. Replacing this structure by the centralized integral controller in~\eqref{eq:capi}, however, notably improves scalability (compare (f) to (e)).   }

\edit{The topologies of the 14 and 118 bus test cases differ from the regular lattices considered in this paper. However, they are planar graphs and can therefore be \emph{embedded} in a 2D lattice. This means that the 2D lattice gives a lower bound on their performance, see~\cite{Tegling2017b} for an elaboration. We can therefore expect, at best, a linear scaling of the normalized \hn norm with noisy distributed integral control, but logarithmic with centralized integral control.      }



\section{Discussion}
\label{sec:discussion}
We conclude this paper by a brief discussion of our findings and some directions for future work.
\subsection{Control design for finite size networks}
The focus of this paper has been to characterize limitations of distributed integral control in terms of the scaling of \hn performance to large networks. We showed that such limitations arise due to noisy measurements, and can only be alleviated by asymptotically increasing the amount of inter-nodal alignment between controllers. 

Naturally, any real-world application will have a finite number of nodes, and the controller can thus always be tuned for acceptable performance. Our results imply, however, that such a tuning cannot be done independently of the network size. Therefore, even though the controller is implemented in a distributed fashion, its tuning requires global knowledge.

\subsection{Alternative dynamic feedback {and nodal dynamics}}
A natural question to ask is whether there are other distributed integral control designs for the system~\eqref{eq:system}, which are not subject to the same performance limitations. In general, one can set   
$u = z$,  $\dot{z} = Az + Bx + Cv^m,$ 
where $A,B,C$ are feedback operators. Even with a controller on this form it is possible to show 
a result in line with Proposition~\ref{cor:bestscaling}.
It is important to note, however, that the important property of \emph{load sharing} is lost if $B \neq 0$. Load sharing implies that the control effort in regulating a constant disturbance is shared equally (or proportionally) between controllers. Therefore, even though a controller with $B\neq 0$ can achieve the same performance as static feedback,  it may not be desirable in practice.  

\edit{Another important question is that of higher-order nodal dynamics (here they were of second order). In this case, the situation differs depending on whether there is absolute feedback from additional states (such as damping terms) or only relative feedback. Relative localized feedback, as we considered here, does not scale well to large networks with high-order dynamics -- stability tends to be lost~\cite[Chapter 5]{TeglingThesis}. It is, however, an open question whether more intricate nodal dynamics (with further absolute states) could be exploited to alleviate performance limitations.}


\subsection{Network heterogeneity}
The results in this paper were derived under the assumption of spatial invariance. This allowed us to obtain closed-form expressions for~\hn norms and their scalings. \edit{In particular, the assumption enabled an analysis of a type of multiplex controller structure (see Figure~1), where the communication network can differ in connectivity from the underlying physical network. }

\edit{Most network applications, however, have other topological structures.} In such more general network topologies, it is often possible to derive performance bounds on similar forms as the results herein using embedding arguments. \edit{See~\cite{Tegling2017b} for a more elaborate discussion. This means that performance limitations that apply to lattices also apply to general networks that can be embedded in them. This was also demonstrated through examples in Section~\ref{sec:examples}.}
\edit{ Otherwise, the correct generalization of the notion of spatial dimension, which is important for the scalings discussed here, was addressed in~\cite{Patterson2014} but remains an open research question. }


\section*{Acknowledgements}
\vspace{-2mm}
We wish to thank Bassam Bamieh and Hendrik Flamme for many insightful discussions related to this work. Funding support from the Swedish Research Council under Grants 2013-5523 and 2016-00861 is also acknowledged. 



\section*{Appendix}
\vspace{-2mm}
\subsection*{{ Proof of Lemma~\ref{lem:q}}}
\edit{Recall that $\hat{a}(\theta) = -\sum_{k \in \mathbb{Z}^d}a_k(1-\cos(\theta\cdot k)) $ and that the smallest wavenumber is $\theta_{\min} = {2\pi}/{L}$.
We establish a bound for $\hat{a}(\theta = {2\pi}/{L})$ as follows (dropping the subscript of~$q_A$): 
\begin{small}
\begin{align}
\nonumber
&\sum_{k =-{q}}^{q} \!\! a_k \! \left( \! 1\!-\! \cos\left(\!\frac{2\pi k}{L} \!\right)\!\!\right)\\ \nonumber
= & ~(a_1+a_{-1})\!\left(\!1-\!\cos\left(\frac{2\pi }{L} \right)\!\!\right) +  \ldots 
+ (a_{q}+a_{-q})\!\left(\!1-\!\cos\left(\frac{2\pi q}{L} \right)\!\!\right) \\ \nonumber
\ge &~
2a_{\min}\frac{2}{\pi^2}\!\!\left(\!\! \left(\frac{2\pi}{L} \right)^2 \!+ \!\left( \! \frac{2\pi \!\cdot\! 2}{L}\! \right)^2 \!+  \ldots + \left(\frac{2\pi q}{L} \right)^2  \right)\\ =&~ \frac{16a_{\min}}{L^2}\!\left( 1^2\! + 2^2 \!+ \! \ldots \!+ q^2 \right)  = \frac{16 a_{\min}}{L^2}\frac{q(q+1)(2q+1)}{6}, \label{eq:qbound}
\end{align}
\end{small}
where the first inequality follows from the fact that $1-\cos x \ge 2/\pi^2 x^2$ if $x \in [-\pi,\pi]$ and the last equality from the expression for a sum of a sequence of squares. Now, if $q = \bar{c} L^{2/3}$, where $\bar{c}>0$ is a fixed constant, then \eqref{eq:qbound} is lower bounded by $\frac{16 a_{\min}}{L^2}\frac{2\bar{c}^3L^2}{6} = \frac{32\bar{c}^3a_{min}}{6}=:\delta$, which is a positive constant independent of $L$. The lemma follows. }

\section*{References}
  \bibliographystyle{elsarticle-num} 
  \bibliography{emmasbib2015,BassamBib}





\end{document}